\documentclass[letter]{IEEEtran}

\usepackage{ifpdf}

\ifCLASSINFOpdf
   \usepackage[pdftex]{graphicx}

\else

\fi

\usepackage[cmex10]{amsmath}

\usepackage{times,amsmath,epsfig,latexsym,multirow,array,amssymb,multicol,stfloats,textcomp}
\usepackage{times,amsmath,epsfig,latexsym,multirow,array,amssymb,multicol,stfloats,textcomp,color,colortbl}
\usepackage{algpseudocode}
\usepackage{algorithm}
\usepackage[applemac]{inputenc}
\usepackage{float}
\usepackage{placeins}
\usepackage[table]{xcolor}

\usepackage{nomencl}
\makenomenclature
\usepackage{eurosym}

\begin{document}

\title{
A Comprehensive Methodology for Volt-VAR Optimization in Active Smart Grids
}

\author{Mohammad~Ghasemi~Damavandi,~\IEEEmembership{Student Member,~IEEE},
Vikram~Krishnamurthy,~\IEEEmembership{Fellow,~IEEE}, and
{Jos\'e}~R.~Mart\'i,~\IEEEmembership{Life Fellow,~IEEE}
\thanks{The authors are with the Department of Electrical and Computer Engineering, University of British Columbia (UBC), Vancouver, BC V6T 1Z4, Canada (e-mail: \{mghasemi, vikramk,  jrms\}@ece.ubc.ca).}}

\maketitle

\IEEEpeerreviewmaketitle

\begin{abstract}
This paper considers the problem of Volt-VAR Optimization (VVO) in active smart grids. Active smart grids are equipped with distributed generators, distributed storage systems, and tie-line switches that allow for topological reconfiguration. In this paper, the joint operation of the remotely controllable switches, storage units, under load tap changers, and shunt capacitors is formulated for a day-ahead operation scheme. The proposed VVO problem aims at minimizing the expected active power loss in the system and is formulated as a mixed-integer quadratic program. The stochasticity of the wind power generation is addressed through a first-order Markov chain model. Numerical results on a $33$-node, $12.66$ kV active smart grid using real data from smart meters and wind turbines are presented and compared for various test cases.

\end{abstract}




\section{Introduction}

\IEEEPARstart{A}{ctive} smart grids are equipped with renewable energy sources distributed in the system. In spite of various advantages of these non-dispatchable Distributed Generators (DGs), they can cause dramatic changes in the voltage profile of the system due to high variability and intermittence. In addition to DGs, active smart grids may also benefit from Distributed Storage Systems (DSS) which can cope with the inherent intermittence of renewable energy sources. By using DSS, the Smart Grid Operator (SGO) can also buy the excess energy from DG owners during off-peak hours and sell it back to the network when the demand is higher. If the price of energy is dependent on the total load of the system, e.g. under time of use pricing schemes, this will result in an arbitrage gain for the SGO. It will also be beneficial to the DG owner as it can prevent the curtailment of the non-dispatchable energy sources \cite{High}. If the system undergoes a fault, the SGO can also employ these DSS during system restoration \cite{restoration} to minimize the energy-not-supplied. The DSS can also contribute to load shaping \cite{LoadShaping}, active loss reduction, voltage control \cite{PrimaryVoltage}, \cite{IntegratingElectrical}, ancillary services, and power smoothing for solar arrays \cite{TheRole}.

In addition to DGs and DSS, active smart grids are equipped with remotely controllable, tie-line switches that allow for topological reconfiguration \cite{RobustMeter}. Use of remotely controllable switches enables the SGO to optimize the configuration of the system even in daily or hourly intervals \cite{Valueof}-\nocite{SegmentedTime, MHCjournal, SwitchingFrequency}\cite{TimeVarying}. In practice, network reconfiguration can be exploited for a variety of objectives including load balancing \cite{NetworkReconfiguration}, reliability improvement \cite{reliability}, service restoration \cite{restoration_reconfig}, voltage profile improvement and power loss minimization \cite{ConvexModels}, \cite{PowerLoss}.
In systems with high penetration of non-dispatchable sources, the reconfiguration routine can also be employed to mitigate the effect of intermittent renewable generation. 

Due to the high variability of non-dispatchable sources, the Volt-VAR Optimization (VVO) problem will be a very important issue in active smart grids. The VVO problem, in the classical sense, is the optimal operation of the capacitor banks and Under Load Tap Changers (ULTCs). The objective of the VVO techniques is to optimize some metric of the system performance such as active power loss or voltage profile. The VVO problem in modern smart grids, however, should also take into account the operation of the DSS and remotely controllable switches as additional degrees of freedom. 

Several works in the literature have considered the VVO problem for distribution systems \cite{HamedAhmadi}-\nocite{PeakLoadRelief, SSandVVO, SensitivityBased, MPCbased }\cite{IntegrationDGVVO}. In \cite{VVO1}, a VVO technique using oriented discrete coordinate descent method for minimization of power loss or the number of control actions is proposed. In \cite{ThroughDG}, a VVO methodology is presented which minimizes the total reactive power supplied by the DGs. In \cite{VVONiknam}, a scenario-based multiobjective method for daily VVO is presented which considers renewable energy sources. Voltage control and VAR optimization methodologies are also presented in \cite{Acomprehensive} and \cite{MISOCP_VAR} that consider feeder reconfiguration in addition to classical VVO devices. However, these papers do not consider the DSS in their VVO methodologies. 

This paper presents a comprehensive formulation of the VVO problem in active smart grids considering wind turbines, DSS, capacitor banks, ULTCs, and feeder reconfiguration. The VVO problem is formulated as a mixed-integer, quadratic program to be solved using branch-and-bound methods. Unlike \cite{Acomprehensive}-\cite{MISOCP_VAR}, the objective of the presented methodology is to minimize the \emph{expected} active power loss in the system. In order to formulate the expected power loss of the system, the stochasticity of the wind power generation is addressed using a first-order Markov Chain (MC) model \cite{MCMC}. Simulation results on a $33$-node, $12.66$ kV reconfigurable smart grid using real data from wind turbines and smart meters are presented and discussed. Various test cases are considered in the simulations to compare the impact of different VVO equipment on the system loss and voltage profile.

The rest of the paper is organized as follows. Sec. \ref{Preliminaries} provides some preliminaries about the power flow equations in radial distribution systems as well as an MC model for wind power generation. Sec. \ref{Fixed} presents the formulation of the VVO problem considering DSS, capacitor banks, ULTCs, and network reconfiguration. Sec. \ref{Simulations} provides the numerical results on a reconfigurable test system and by using real data from smart meters and wind turbines. Finally, Sec. \ref{Conclusion} concludes the paper.


\section{Preliminaries} \label{Preliminaries}

This section reviews the DistFlow equations of radial distribution systems \cite{NetworkReconfiguration}, \cite{CapacitorPlacement} and extends them to accommodate the DGs and bi-directional flow of power. Also, an MC model for wind power generation based on \cite{MCMC} is reviewed.

Consider a radial distribution system with $N$ nodes and let $v_n = |v_n| e^{j\theta_n}, n=1,2,\dots,N$ be the voltage phasor of the $n$th node, where $j = \sqrt{-1}$. Let $s_n = p_n + j q_n$ denote the complex power demand of a generic node $n$. Similarly, let $s^G_n = p^G_n + j q^G_n$ be the complex power generation of node $n$. For any radial configuration of the system, let $\mathcal{L}$ be the set of all lines which connect two nodes of the system together. Note that $\mathcal{L}$ uniquely specifies the configuration of the system and, therefore, we denote the configuration of the system by $\mathcal{L}$. With some abuse of notation, we also use $(n,m) \in \mathcal{L}$ to specify that node $n$ is connected to node $m$ in configuration $\mathcal{L}$. In addition, let $\mathcal{N}^\mathcal{L}_n$ be the set of all nodes which are connected to node $n$ via a line in $\mathcal{L}$. Also, let $\mathcal{S}_\text{sub}$ be the set of substations in the system. Throughout the paper, it is assumed that the renewable sources are controlled to have a unit power factor. Likewise, it is assumed that the DSS only store and supply active power, and reactive power demands will be compensated by capacitors as needed. 

\subsection{DistFlow Equations} 

Suppose that node $n$ is connected to node $m$ in a radial configuration via line $l \in \mathcal{L}$ with impedance $z_l = r_l + j x_l$. For any instant of time, let $s_{nm} = p_{nm} + j q_{nm}$ be the complex power that flows from node $n$ towards node $m$. Also, let $s_{mn} = p_{mn} + j q_{mn}$ denote the complex power that flows from node $m$ towards node $n$. Then we have:
\begin{align}
s_{nm} + s_{mn} &= z_l  \frac{|s_{nm}|^2}{|v_n|^2} \nonumber \\
&= z_l  \frac{|s_{mn}|^2}{|v_m|^2} \label{1} \\
v_m &= v_n - z_l \frac{s_{nm}^*}{v_n^*} \label{2}.
\end{align}
Separating the active and reactive parts in \eqref{1} and the voltage magnitudes in \eqref{2}, yields:
\begin{align}
p_{nm} + p_{mn} &= r_l  \frac{p_{nm}^2+q_{nm}^2}{|v_n|^2} \nonumber \\
& = r_l  \frac{p_{mn}^2+q_{mn}^2}{|v_m|^2} \label{3}
\end{align}
\begin{align}
q_{nm} + q_{mn} &= x_l  \frac{p_{nm}^2+q_{nm}^2}{|v_n|^2} \nonumber \\
& = x_l  \frac{p_{mn}^2+q_{mn}^2}{|v_m|^2} \label{4}
\end{align}
\begin{align}
|v_m|^2 &= |v_n|^2 - 2 (r_l p_{nm} + x_l q_{nm}) \nonumber \\
&~~~~~~+ (r_l^2+x_l^2) \frac{p_{nm}^2+q_{nm}^2}{|v_n|^2} \label{5}.
\end{align}
Equations \eqref{3}-\eqref{5} together with the following generation-demand equations
\begin{align}
p^G_n - p_n &= \sum_{m \in \mathcal{N}^\mathcal{L}_n} p_{nm} \label{6} \\
q^G_n - q_n &= \sum_{m \in \mathcal{N}^\mathcal{L}_n} q_{nm} \label{7}
\end{align}
are called the \emph{DistFlow} equations and provide a full AC power flow model for a radial distribution system.
Note that these equations are slightly different from the DistFlow equations in \cite{ConvexModels} and \cite{NetworkReconfiguration} in that they support bi-directional 
flow of power which is required for radial systems with distributed generations. Accordingly, the active and reactive powers that flow from one node to another can take both positive and negative values.

\subsection{Markov Chain Model of Wind Power Generation}

The wind power generation can be modeled as a first-order homogeneous Markov Chain (MC) as proposed in \cite{MCMC}. To that end, the range of the power generation by the wind turbine is divided into $S$ intervals, each denoted by a state level $x_i, i=1,2, \dots, S$. Here, $x_i$ is the mean power generated by the wind turbine when the power generation is in the $i^{\text{th}}$ interval. Next, a transition probability $0 \leq \theta_{ij} \leq 1$ is assigned to the event that the wind power generation goes from state $i$ at hour $h$ to the state $j$ at hour $h+1$. These probabilities form the transition probability matrix $\mathbf{T} = [\theta_{ij}]$. The matrix $\mathbf{T}$ is a stochastic matrix satisfying $\mathbf{T} \mathbf{1} = \mathbf{1}$, where $\mathbf{1}$  denotes a vector with all elements equal to one.

Let $\boldsymbol{\pi}(h) = [\pi_1(h), \pi_2(h), \dots, \pi_S(h)]'$ be the state probability vector at hour $h$. Then, based on the Chapman-Kolmogorov equation, the probability of being in each state at hour $h+\tau$ will be $\boldsymbol{\pi}(h+\tau) = \left(\mathbf{T'}\right)^{\tau} \boldsymbol{\pi}(h)$. 
In practice, the realization of the wind power generation is observed at the beginning of the optimization process. Therefore, all the elements of $\boldsymbol{\pi}(0)$ but one are equal to zero. The initial state probability vector $\boldsymbol{\pi}(0)$ is used to obtain subsequent state probability vectors.

The entries of the matrix $\mathbf{T}$ can be estimated from real data using a Maximum-Likelihood (ML) approach. Let $n_{ij}$ be the number of times that the MC goes from state $i$ to state $j$ in the available dataset. Then, the ML estimate of $\theta_{ij}$ is computed to be:
\begin{equation}
\hat{\theta}_{ij} = \frac{n_{ij}}{\sum_j n_{ij}}.
\end{equation}

Unlike wind power generation, the loads of the system usually have daily patterns that can be forecast with acceptable accuracy based on season, day of the week, time of the day, and weather conditions. Therefore, similarly to \cite{Relief}\nocite{HamedAhmadi,Integration,InverterLess}-\cite{TimeInterval}, typical load patterns obtained using smart meter measurements are employed for day-ahead VVO in this research. In particular, the load pattern of each node is simply obtained by computing the average load of the node over a history of a prescribed length, e.g., one month. Once the day-ahead schedule of the system equipments are obtained based on typical load patterns, finer adjustments can be made in real time based on the actual measurements of the loads \cite{HamedAhmadi, HybridVVO}.
 
\section{Formulation of the Volt-VAR Optimization Problem}\label{Fixed}

This section formulates the VVO problem in active smart grids considering wind generation, DSS, ULTCs, capacitors, and feeder reconfiguration. The terms DSS and \emph{storage units} will be used interchangeably in this paper.

\subsection{The Objective Function of the VVO Problem}

In this paper, the objective function of the VVO problem is considered to be the expected active power loss in the system. Based on the DistFlow equations, the total active power loss in the system at hour $h$ is given by:
\begin{equation}
\text{Loss}(h) = \sum_{(n,m) \in \mathcal{L}} r_l  \frac{p_{nm}^2(h)+q_{nm}^2(h)}{|v_n(h)|^2}.
\end{equation}
Approximating the nodal voltages by $1$ p.u. in the right-hand side of \eqref{3}, the total active power loss in the system at hour $h$ will be given by the following quadratic form:
\begin{equation}\label{LossAppr}
\text{Loss}(h) \approx \sum_{(n,m) \in \mathcal{L}} r_l  \left(p_{nm}^2(h)+q_{nm}^2(h) \right).
\end{equation} 
The active power loss in the system at each hour is a random variable due to the stochasticity of the wind power generation. Therefore, the VVO problem should try to minimize the expected active power loss in the system as:
\begin{equation}
\min \mathbb{E} \left\{\frac{1}{H} \sum_{h=1}^H \text{Loss}(h) \right\}, \label{Loss}
\end{equation}
where $\mathbb{E}\{ \cdot \}$ is the expectation operator. Also, $H$ is the horizon of the VVO problem which, in the case of a day-ahead approach, equals $24$. One can employ the MC model of the wind power generation to compute the expectation in \eqref{Loss} as follows:
\begin{align}\label{expand}
&\mathbb{E} \left\{\frac{1}{H} \sum_{h=1}^H \text{Loss}(h) \right\} = \frac{1}{H} \sum_{i=1}^S \sum_{h=1}^H  \pi_i(h) ~ \text{Loss}(h;i),
\end{align}
where $\text{Loss}(h;i)$ is the total power loss at hour $h$ assuming that the wind power generation is in the $i^{\text{th}}$ state. It is seen from \eqref{LossAppr} and \eqref{expand} that the objective function of the VVO problem is quadratic and convex.

\subsection{Power Flow Equations, Distributed Storage Systems, Capacitors, and ULTCs}

The VVO problem needs to be solved subject to a series of constraints. Below, the constraints corresponding to power flow equations as well as the operation of DSS, capacitors, and ULTCs are formulated.

A full set of AC power flow equations can be described by \eqref{3}-\eqref{5} along with the generation-demand equations. To include the charging and discharging strategies of the storage units as well as the reactive power injection of  capacitors, the generation-demand equations \eqref{6}-\eqref{7} have to be modified as follows:
\begin{align}
&p^G_n(h) - p_n(h) - p^\text{DSS}_n(h) = \sum_{m \in \mathcal{N}^\mathcal{L}_n} p_{nm}(h) , \nonumber \\
&\qquad \qquad \qquad \qquad \qquad \qquad \qquad ~~~h \in \mathcal{H}_o, ~\forall n, \label{15} 
\end{align}
\begin{align}
&p^G_n(h) - p_n(h) + p^\text{DSS}_n(h) = \sum_{m \in \mathcal{N}^\mathcal{L}_n} p_{nm}(h), \nonumber \\
&\qquad \qquad \qquad \qquad \qquad \qquad \qquad ~~~h \in \mathcal{H}_p,  ~\forall n, \label{14} 
\end{align}
\begin{align}
&q^G_n(h) - q_n(h) = \sum_{m \in \mathcal{N}^\mathcal{L}_n} q_{nm}(h), \nonumber \\
&\qquad \qquad \qquad \qquad \qquad \qquad ~~h \in \mathcal{H}_o \cup \mathcal{H}_p,  ~\forall n. \label{20}
\end{align}
In this formulation, $p^\text{DSS}_n(h)$ is the amount of average power that the storage unit installed on node $n$ stores from or supplies to the grid at hour $h$. Also, $q^G_n(h)$ is the amount of reactive power generation of the capacitor installed on node $n$ at hour $h$. 
Moreover, $\mathcal{H}_p$ and $\mathcal{H}_o$ denote the set of peak and off-peak hours, respectively.

The constraints corresponding to the rated power limit of the storage units during charging and discharging periods can be written as:
\begin{align}
\mathbf{0}_N \leq \mathbf{p}^\text{DSS}(h) \leq \mathbf{p}^{\text{DSS},max}, ~~~~\forall h, \label{17}
\end{align}  
where $\mathbf{0}_N$ is an all-zero vector of length $N$. Also,
$
\mathbf{p}^\text{DSS}(h) = [p^\text{DSS}_1(h), p^\text{DSS}_2(h), \dots, p^\text{DSS}_N(h)]^T
$
is the vector of power storage and injections of the storage units at hour $h$, and
$
\mathbf{p}^{\text{DSS},max} = [p^{\text{DSS},max}_1, p^{\text{DSS},max}_2, \dots, p^{\text{DSS},max}_N]^T
$
is the vector of power ratings of the storage units. If a node $n$ is not equipped with a storage unit, then its power rating is set to zero, that is:
\begin{equation}
p^{\text{DSS},max}_n = 0, ~~~ n \notin \mathcal{S}_\text{DSS},
\end{equation}
where $\mathcal{S}_\text{DSS}$ is the set of nodes equipped with DSS.

The constraints corresponding to the capacity of the storage units can be formulated as:
\begin{align}
\sum_{h \in \mathcal{H}_o} \mathbf{p}^\text{DSS}(h) &= \frac{\gamma_{\text{DOD}}}{\beta_{ch}} \mathbf{b},  \label{18} \\
\sum_{h \in \mathcal{H}_p} \mathbf{p}^\text{DSS}(h) &= \beta_{dis} \gamma_{\text{DOD}}  \mathbf{b}, \label{19} 
\end{align}
where $\mathbf{b} = [b_1,b_2, \dots, b_n]^T$ in kWh is the vector of DSS capacities installed in the system. If a node $n$ is not equipped with DSS, then $b_n = 0$. Also, $\gamma_{\text{DOD}}$ is the Depth of Discharge (DOD) of the storage units and $0 < \beta_{ch} < 1$ and $0 < \beta_{dis} < 1$ are the charging and discharging efficiencies of the DSS technology, respectively. Hence, $\beta_{ch} \beta_{dis}$ is the round-trip efficiency of the storage units. 
Note that the first time that a storage unit of size $b_n$ is connected to node $n$, a total energy of $\frac{b_n}{\beta_{ch}}$ kWh will be absorbed to charge it up. After the first charging period, always $(1-\gamma_{\text{DOD}}) b_n$ kWh of energy remains in the unit while the remaining $\gamma_{\text{DOD}} b_n$ kWh is absorbed during charging periods. This explains the factor $\gamma_{\text{DOD}}$ on the right-hand side of \eqref{18} and \eqref{19}. At off-peak hours when the storage unit gets charged, the total energy that is fed to the unit times the charging efficiency coefficient of the unit has to be equal to $\gamma_{\text{DOD}} b_n$, hence \eqref{18}. At peak hours when the storage unit gets discharged, the total dischargable energy of the unit times its discharging efficiency coefficient has to be equal to the total energy that is injected to the grid, hence \eqref{19}.
Note that our formulation assumes that the storage units are charged during off-peak hours and discharged during peak hours. That is, the charging and discharging pattern of the storage units is only once a day. This assumption ensures a longer lifetime for the storage units and has been made extensively in the literature \cite{GridSupport}, \cite{FlexibleOptimal}. 

The shunt capacitors can inject reactive power into the feeder. This reactive power is modeled as:
\begin{align}
q^G_n(h) &= c_n(h) q^{cap}_n, ~~\quad\quad \qquad n \in \mathcal{S}_\text{cap},\\
c_n(h) &\in \{0,1,2, \dots,  c_n^{max}\},~~~~ \forall h. \label{22}
\end{align} 
where $q^{cap}_n$ is the reactive power injection by a unit module of the capacitor bank installed on node $n$ and $c_n(h)$ is the number of modules connected at hour $h$. Also, $\mathcal{S}_\text{cap}$ is the set of nodes equipped with a shunt capacitor and $c_n^{max}$ is the maximum number of modules available in the capacitor bank installed on node $n$. 

A ULTC that is installed on the branch connecting node $n$ to node $m$ with branch impedance of $z_l$ can be modeled as follows \cite{EnergyStorageDevices, MISOCP_VAR}. An ideal transformer is modeled between node $n$ and an auxiliary node $m'$, in series with the impedance $z_l$ between node $m'$ and node $m$. Therefore, the ULTC operation can be formulated as:
\begin{align}
&|v_{m}|^2 = |v_{m'}|^2 - 2 (r_l p_{{m'}m} + x_l q_{{m'}m}) \nonumber \\
&\qquad \qquad+ (r_l^2+x_l^2) \frac{p_{{m'}m}^2+q_{{m'}m}^2}{|v_{m'}|^2} .  \\
&|v_{m'}(h)| = (1 + \gamma \Delta) |v_n(h)|, \label{ULTC}
\end{align}
where $\gamma$ is the ULTC tap position and $\Delta$ is the tap step size. For instance, if the transformer ratio is between $1-a$ and $1 + a$ with step size $\Delta$, then
\begin{equation}\label{tap}
\gamma \in \left\{\gamma_i = i \big| i = -\frac{a}{\Delta}, -\frac{a}{\Delta}+1, \dots, \frac{a}{\Delta}-1, \frac{a}{\Delta} \right\}.
\end{equation}
For the sake of simplicity, no explicit index has been used for the ULTC tap position and step size to indicate the corresponding node numbers. 

Finally, the operational constraints on the nodal voltages and feeder ampacities should be formulated. The nodal voltages should be  bounded as:
\begin{align}
\mathbf{v}^{min} \leq |\mathbf{v}(h)| \leq \mathbf{v}^{max}, ~~ ~\forall h \label{9},
\end{align}
where $\mathbf{v}(h) = [|v_1(h)|, |v_2(h)|, \dots, |v_N(h)|]^T$ is the vector of nodal voltage magnitudes at hour $h$ and $\mathbf{v}^{min} = [v_1^{min},v_2^{min}, \dots, v_N^{min}]^T$ and $\mathbf{v}^{max} = [v_1^{max},v_2^{max}, \dots, v_N^{max}]^T$ are the vector of minimum and maximum allowable voltage magnitudes, respectively. Also, the limits on feeder ampacities can be written as:
\begin{align}
p^2_{nm}(h) + q^2_{nm}(h) \leq |s^{max}_{nm}|^2, ~&(m,n) \in \mathcal{L},  ~\forall h \\
p^2_{mn}(h) + q^2_{mn}(h) \leq |s^{max}_{nm}|^2, ~&(n,m) \in \mathcal{L},  ~\forall h \label{10}
\end{align}
where $|s^{max}_{nm}|$ is the maximum apparent power that can flow in the line which connects node $n$ to node $m$.


\subsection{Feeder Reconfiguration}\label{Active}

Modern smart grids are equipped with remotely-controllable tie-line switches for topological reconfiguration. This section formulates the role of feeder reconfiguration on the VVO problem. It is assumed that the reconfiguration routine is performed once in $24$ hours, i.e., in a day-ahead fashion. 

To formulate the constraints corresponding to feeder reconfiguration, the extended formulation of \cite{ConvexModels} is exploited which incorporates the existence of DGs and bi-directional flow of power. Let $\mathcal{L}_\infty$ be the layout of the network, i.e. the configuration of the system when all the tie-line switches are closed. Denote by $|\mathcal{L}_\infty|$ the cardinality of $\mathcal{L}_\infty$. Define $|\mathcal{L}_\infty|$ binary variables $y_{nm} \in \{0,1\}$ associated with each line between any two nodes in the final radial configuration as:
\begin{equation}
y_{nm} = \left\{
\begin{array}{rl}
1 , ~~~~~&\text{if node $n$ is connected to node $m$}  \label{y}\\
0 , ~~~~~&  \text{otherwise}\\
\end{array} \right.
\end{equation}
and stack them in $\mathbf{y}$. Using these notations, the constraints corresponding to the reconfiguration routine can be formulated as follows.
For any $n \in \{1,2, \dots,N \}$, $(n,m) \in \mathcal{L}_\infty$, the generation-demand equations \eqref{15}-\eqref{20} in a reconfigurable distribution system should be augmented with the following constraints:
\begin{align}
-D y_{nm} &\leq p_{nm}(h) \leq D y_{nm}, ~~~~~ \forall h \label{63}\\
-D y_{nm} &\leq q_{nm}(h) \leq D y_{nm}, ~~~~~ \forall h \\
-D y_{nm} &\leq p_{mn}(h) \leq D y_{nm}, ~~~~~ \forall h \\
-D y_{nm} &\leq q_{mn}(h) \leq D y_{nm}, ~~~~~ \forall h \label{64}
\end{align}
where $D \gg 1$ is a large disjunctive constant.

The radiality constraint can be imposed by requiring that the resulting configuration should be loop-free and connected. 
As explained in \cite{Imposing}, for distribution systems with DGs this can be achieved by introducing a set of fictitious loads on the nodes that can be fed only from the substation. Let $k_n$ be the fictitious load on node $n$. Also, let $k_{nm}$ be the fictitious power that flows from node $n$ towards node $m$. Then, the following set of linear constraints are equivalent to the radiality of the system.
\begin{align}
\mathbf{1}^T \mathbf{y} &= N-1, \label{500}\\
k_n &= \sum_{m \in \mathcal{N}^{\mathcal{L}_\infty}_n} k_{nm}, ~~\forall n, (n,m) \in \mathcal{L}_\infty, \label{50} \\
-D &y_{nm} \leq k_{nm} \leq D y_{nm}, \label{51}
\end{align}
where $\mathbf{1}$ is an all-one vector. Note that it is required to impose the radiality constraint for only one hour. Therefore, $k_m$ and $k_{nm}$ do not depend on $h$. 

In addition, it is possible to limit the number of switching actions during the reconfiguration procedure to prevent excessive costs. Let $\bar{\mathbf{y}}$ denote the status of the switches in the current configuration of the system. Also, let $\boldsymbol{\alpha}$ be an auxiliary vector. Then, the following set of linear constraints limits the number of switching actions to $\rho$:
\begin{align}
\mathbf{y} - \bar{\mathbf{y}} &\leq \boldsymbol{\alpha}, \\
\bar{\mathbf{y}} - \mathbf{y} & \leq \boldsymbol{\alpha},  \\
\mathbf{1}^T \boldsymbol{\alpha} &\leq \rho, \label{rho} \\
\boldsymbol{\alpha} &\leq \mathbf{1}. \label{alfa}
\end{align}
Note that for the system to remain connected, opening one switch requires closing another switch. Therefore, $\rho$ should be an even number in \eqref{rho}. 

Finally, if a branch of the system is not equipped with remotely controllable switches, its connection is forced to remain unchanged. That is,
\begin{equation}\label{TieSet}
\qquad \qquad ~~~~ y_{nm} = \bar{y}_{nm}, \qquad (n,m) \notin \mathcal{S}_\text{tie}
\end{equation}
where $ \mathcal{S}_\text{tie}$ is the set of branches equipped with a remotely controllable switch. 

\subsection{Convexification of the VVO Problem}

The VVO problem derived above is non-convex due to some non-linearities in the power flow equations, the multiplication of variables in the right-hand side of \eqref{ULTC}, and the integrality constraints \eqref{22} and \eqref{tap}. However, it is possible to make approximations to the power flow equations \cite{ConvexModels} to come up with a mixed-integer convex optimization problem. To that end, we will need to consider the square of the voltage magnitudes as independent optimization variables. Also, the multiplication of variables in \eqref{ULTC} can be replaced by alternative convex constraints as will be discussed next. 

Although mixed-integer convex optimization problems are NP-hard, the relaxed version of those problems is convex and, hence, they can often be solved in a reasonable time. In practice, mixed-integer convex programs are solved using a combination of a convex optimization technique and an exhaustive search algorithm, such as branch-and-bound methods. Two main characteristics that makes solving mixed-integer convex programs particularly easier are the following. First, when performing the exhaustive search over integer variables along a tree, some branches of the tree can be shown (through solving a relaxed problem) not to include the optimal solution. Therefore, there is no need to follow those branches. Second, if an optimal solution is found before searching all combinations of the integer variables, it may be possible to prove that this solution is optimal through solving a relaxed problem. This is in contrast with non-convex integer programs, where even if the optimal solution is obtained before searching all the combinations, there is no way to prove that this solution is actually optimal. 


To convexify the VVO problem, first, constraint \eqref{ULTC} is replaced by a series of linear constraints using the approach presented in \cite{EnergyStorageDevices}. Since the square of the voltage magnitudes are considered as the optimization variables, \eqref{ULTC} should be written as:
\begin{equation}\label{ULTC2}
|v_{m'}(h)|^2 = (1 + \gamma \Delta)^2 |v_n(h)|^2,
\end{equation}
and
\begin{equation}
(1 + \gamma \Delta)^2 = 1 + 2 \Delta \gamma + \Delta^2 \gamma^2.
\end{equation}
The integer variable $\gamma$ can be expanded using a series of binary variables $\kappa_i \in \{0,1\}$ as:
\begin{align}
&\gamma = \sum_{i=-\frac{a}{\Delta}}^{\frac{a}{\Delta}} \kappa_i \gamma_i, \\
&\sum_{i=-\frac{a}{\Delta}}^{\frac{a}{\Delta}} \kappa_i = 1.
\end{align}
Similarly, 
\begin{equation}
\gamma^2 = \sum_{i=-\frac{a}{\Delta}}^{\frac{a}{\Delta}} \kappa_i \gamma^2_i.
\end{equation}
Therefore, one can rewrite \eqref{ULTC2} as:
\begin{equation}\label{ULTC4}
|v_{m'}(h)|^2 \!=\! |v_n(h)|^2 + 2 \Delta \!\sum_i \! \gamma_i |\tilde{v}_{n,i}(h)|^2 + \Delta^2 \!\sum_i \! \gamma_i^2 |\tilde{v}_{n,i}(h)|^2, 
\end{equation}
where
\begin{equation}\label{multiplication}
|\tilde{v}_{n,i}(h)|^2 = \kappa_i |v_n(h)|^2, ~~~~i = -\frac{a}{\Delta}, \dots, \frac{a}{\Delta}.
\end{equation}
Constraint \eqref{ULTC4} is linear in $|v_m(h)|^2$, $|v_n(h)|^2$, and the $\frac{2 a}{\Delta}+1$ variables $|\tilde{v}_{n,i}(h)|^2$. Constraints \eqref{multiplication} which are multiplications of a binary variable and a bounded continuous variable can each be replaced by the following linear constraints \cite{HamedAhmadi, EnergyStorageDevices}:
\begin{align}
&|v_n(h)|^2 - (1-\kappa_i) (v_2^{max})^2 \leq |\tilde{v}_{n,i}(h)|^2, \label{inja1}\\
&|\tilde{v}_{n,i}(h)|^2 \leq |v_n(h)|^2 - (1-\kappa_i) (v_2^{min})^2, \label{inja2}\\
&\kappa_i (v_2^{min})^2 \leq |\tilde{v}_{n,i}(h)|^2 \leq \kappa_i (v_2^{max})^2, \label{inja3}
\end{align}
If $\kappa_i$ is zero, \eqref{inja3} forces $|\tilde{v}_{n,i}(h)|^2$ to be zero and \eqref{inja1} and \eqref{inja2} become redundant (same as voltage magnitude constraints \eqref{9}). If $\kappa_i$ is one, \eqref{inja1} and \eqref{inja2} force $|\tilde{v}_{n,i}(h)|^2$ to be equal to $|v_n(h)|^2$ and \eqref{inja3} becomes redundant.

To come up with a mixed-integer convex approximation of the VVO problem, it remains to convexify the power flow equations \eqref{3}-\eqref{5}. Note that $r_l (p_{nm}^2+q_{nm}^2)/|v_n|^2$ and $x_l (p_{nm}^2+q_{nm}^2)/|v_n|^2$ are, respectively, the active and reactive power loss on the line $l$ which connects node $n$ to node $m$. Therefore, by ignoring the power loss on the lines, the power flow equations \eqref{3}-\eqref{4} assume the following linear form:
\begin{align}
&p_{nm}(h) + p_{mn}(h) = 0, \label{301} ~~~~~\forall h, (n,m) \in \mathcal{L} \\
&q_{nm}(h) + q_{mn}(h) = 0, \label{401} ~~~~~\forall h, (n,m) \in \mathcal{L}.
\end{align}


Moreover, since $r_l$ and $x_l$ have small values in real distribution systems, one can drop the third term in the right-hand side of \eqref{5}. That is, \eqref{5} can be approximated with the following constraint:
\begin{align}
&|v_m(h)|^2 = |v_n(h)|^2 - 2 \big(r_l p_{nm}(h) + x_l q_{nm}(h) \big), \nonumber \\
&~~ \qquad \qquad \qquad \qquad \qquad \qquad \qquad ~\forall h, (n,m) \in \mathcal{L}. \label{501} 
\end{align}
Notice that if $|v_n|^2$ is considered as an independent variable of the optimization then \eqref{501} is linear.

The VVO problem \eqref{Loss} subject to \eqref{301}-\eqref{501}, \eqref{15}-\eqref{22}, \eqref{9}-\eqref{10}, \eqref{63}-\eqref{TieSet}, \eqref{ULTC4}, and \eqref{inja1}-\eqref{inja3} is a mixed-integer quadratic problem and optimizes the joint operation of the DSS, shunt capacitors, ULTCs, and remotely controllable reconfiguration switches for expected power loss minimization.

\section{Numerical Results}\label{Simulations}

In this section, the presented VVO methodology is tested on an active distribution system and several cases with different VVO equipment availability are compared. 

\subsection{The Setting of the Simulations}

The VVO problem for loss minimization is considered for a $33$-node, $12.66$ kV \cite{node_32} active distribution system. In the presented simulations, the Advanced Metering Infrastructure (AMI) data released by the Commission for Energy Regulation (CER) \cite{CER} is utilized to model the loads. The dataset provided by CER is from the Electricity Customer Behaviour Trail study and has been collected from $5000$ smart meters in Ireland from July $14$, $2009$ to December $31$, $2010$. This dataset was received by authors from Irish Social Science Data Archive (ISSDA) \cite{ISSDA}. The reactive power demands are modelled by assuming a constant power factor for the nodes. The horizon of the VVO problem is considered to be $24$ hours corresponding to a day-ahead scenario.

A wind turbine of a rated power of $1000$ kW is installed on node $15$. The time series of wind generation is obtained from \cite{WindData}. The rated capacity of the wind turbine is divided into $10$ states for MC modeling. Based on the real data of the wind turbine, the mean power generation in the states (i.e., the state levels) is $53.96$ kW, $147.16$ kW, $246.9$ kW, $347.47$ kW, $446.96$ kW, $549.96$ kW, $648.19$ kW, $745.16$ kW, $848.05$ kW, $948.07$ kW, respectively.

Two storage units of capacities $200$ kWh and $300$ kWh are installed on the nodes $14$ and $15$, respectively. The power rating of both units is assumed to be $100$ kW. The specifications of the DSS are obtained from \cite{BES} and \cite{Integration} and are as follows. The charging and discharging efficiency of the storage units is $\beta_{ch} = \beta_{dis} = 0.85$. The DOD of the storage units is assumed to be $\gamma_\text{DOD} = 0.75$. 

Two shunt capacitors are installed on the nodes $11$ and $25$. The reactive power rating of the the capacitors are assumed to be $400$ kvar consisting of $4$ modules of size $100$ kvar. A ULTC with the ratio in the range of $0.9 \sim 1.1$ and with tap step sizes of $0.01$ is installed between nodes $6$ and $26$. The branches equipped with remotely controllable tie line switches are the normally-open branches $(8,21)$, $(9,15)$, $(12,22)$, $(18,33)$, $(25,29)$ as well as branches $(6,7)$, $(10,11)$, $(14,15)$, $(29,30)$. The maximum number of switching actions per day is considered to be $\rho = 6$ which is equivalent to opening $3$ switches and closing $3$ other switches. The minimum and maximum voltage limits for all nodes is considered to be $0.94$ p.u. and $1.06$ p.u., respectively. The disjunctive parameter considered for the reconfiguration routine is equal to $100$.


The VVO technique presented in the paper is conducted on the system under study for $N_{\text{day}} = 30$ days and the average results are reported. Typical load patterns obtained by averaging over a history of one month are used for optimization. Also, the MC model of the wind power generation is used during the optimization. Once the optimized operation of the VVO equipment is obtained in a day-ahead fashion, the system is tested with actual loads and wind power generation to compute the nodal voltages and the active power loss. In addition to the average active power loss of the system, four more metrics of quality are also reported. These metrics are the average peak load of the system, the average minimum voltage in the feeder, the average maximum voltage in the feeder, and the average voltage spread in the feeder. The average voltage spread in the feeder is computed to measure the voltage deviation in the system and is defined as:
\begin{equation}
\text{Voltage Spread} = \frac{1}{N_{\text{day}}H} \sum_{i=1}^{N_{\text{day}}} \sum_{h=1}^{H} 
\left[v_\text{max}(h) - v_\text{min}(h) \right],
\end{equation}
where
\begin{align}
v_\text{max}(h) &= \max_n |v_n(h)|, \\
v_\text{min}(h) &= \min_n |v_n(h)|.
\end{align}

All simulations are done using MATPOWER \cite{MATPOWER} and the optimization problems are solved using the package CVX: Software for Disciplined Convex Programming, version 2.1 \cite{CVX} bundled with Mosek, version 7.1. 

\subsection{Results}

As explained in Sec. \ref{Preliminaries}, a linearized version of DistFlow equations is employed in this paper which allows for bi-directional flow of power in a radial system. To first examine the accuracy of the linearized DistFlow equations, a study is conducted on the system. Unlike the VVO optimization procedure, it is assumed in this study that the actual load and wind power generation is available for solving the power flow equations. This study in conducted in the default configuration of the system without shunt capacitors and DSS. The position of the ULTC is fixed at a ratio of $1.05$. The simulations are done for $30$ days, $24$ hour each, and the results are compared with that of the Newton's AC power flow method. The results show that the average difference in the voltage magnitudes between the linearized extended DistFlow equations and the Newton's method is $3.35 \times 10^{-4}$ p.u.. Moreover, the average maximum nodal error of the linearized extended DistFlow equations turns out to be $5.36 \times 10^{-4}$ p.u.. Also, the average error and the average maximum error under peak load are $5.39 \times 10^{-4}$ p.u. and $8.55 \times 10^{-4}$ p.u., respectively. In addition, the average error in the active power loss computed using the linearized DistFlow equations turns out to be $1.56$ kW or $3.3 \%$. Fig. 1 depicts the voltage profile of the system for a typical day and for a peak and an off-peak hour. The figure contrasts the voltage profile obtained by the linearized extended DistFlow equations and that of the Newton's method.

\begin{figure}[btp!] 
\center
\includegraphics[height=2.6in, width = 3.4in]{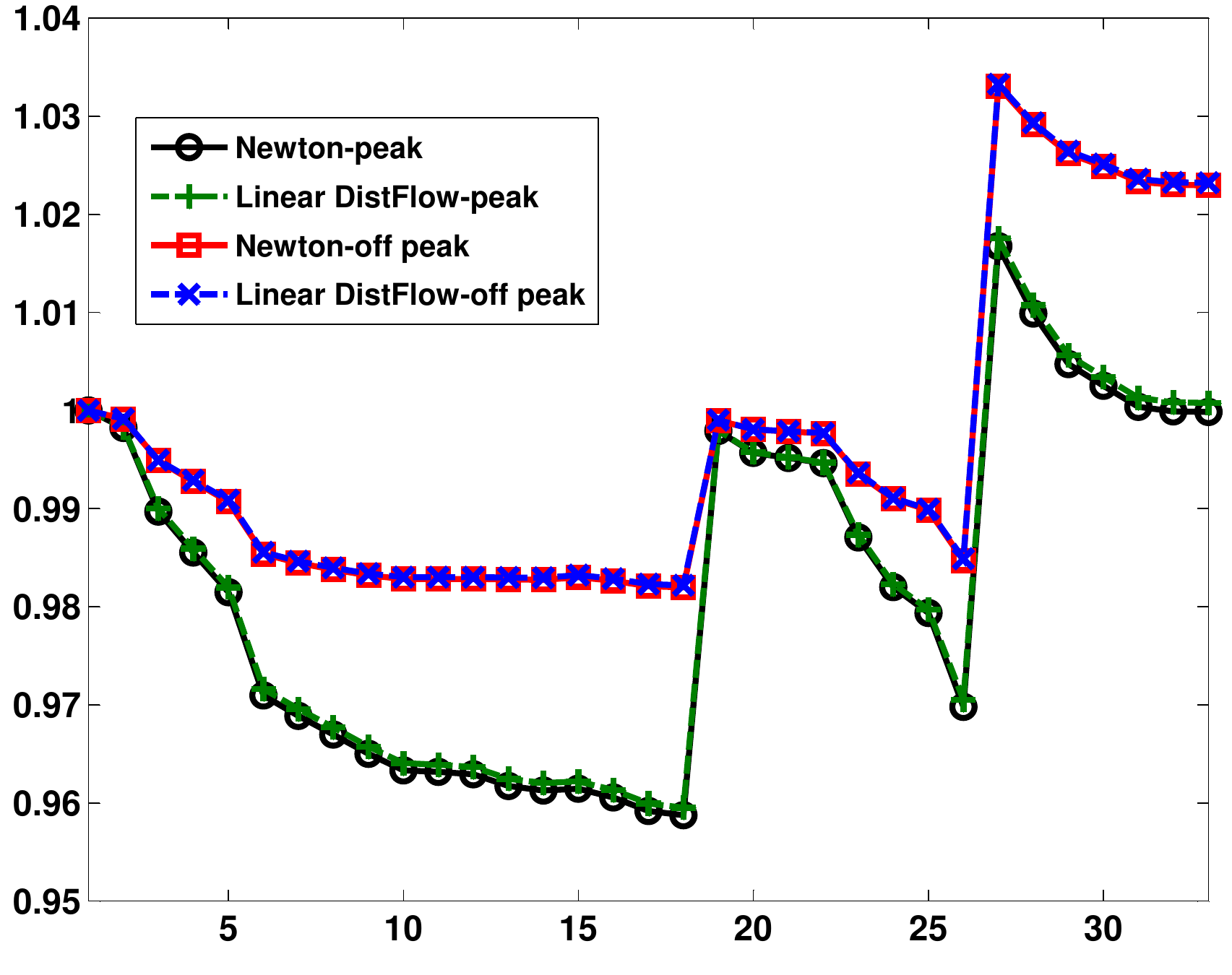}
\caption{Voltage profile of the system on a typical day at a peak and an off-peak hour. Newton's AC power flow versus linearized extended DistFlow equations.}
\end{figure}

\begin{table*}[tbp!]
\caption{Performance of the Optimal VVO Solution Under Different Settings}
\centering
  \begin{tabular}{c cc cccccccccc}
   \noalign{\hrule height 1pt} 
Setting  & Default &\hspace{-1.2mm}ULTC &\hspace{-1.5mm}Reconfig.  &Capacitor 
&\begin{tabular}{c}Capacitor \\DSS \end{tabular}
&\hspace{-1.2mm} \begin{tabular}{c}ULTC \\ Reconfig. \end{tabular} 
&\hspace{-1.2mm} \begin{tabular}{c}ULTC \\ Reconfig. \\DSS \end{tabular}
&\hspace{-1.2mm} \begin{tabular}{c}ULTC \\ Reconfig. \\Capacitor \end{tabular}
& \hspace{-1.2mm}\begin{tabular}{c}ULTC \\ Reconfig. \\Capacitor\\DSS \end{tabular}\\

\noalign{\hrule height 1pt} 
\multirow{2}{*}{} 
   Loss (kW) &$37.94$ &$37.29$ &$31.30$ &$30.47$ &$30.15$  &$30.66$&$30.55$ &$26.97$ &$25.69$ \\
\cellcolor{gray!25}Reduced by (\%)  &\cellcolor{gray!25} -  &\cellcolor{gray!25}$1.72$ &\cellcolor{gray!25}$17.51$ &\cellcolor{gray!25}$19.69$  &\cellcolor{gray!25}$20.52$  &\cellcolor{gray!25}$19.19$ &\cellcolor{gray!25}$19.47$ &\cellcolor{gray!25}$28.91$ &\cellcolor{gray!25}$32.27$\\

\noalign{\hrule height .1pt}
\multirow{2}{*}{} 
   Peak Load (MW) &$2.256$ &$2.255$ &$2.241$ &$2.241$ &$2.169$  &$2.240$&$2.174$ &$2.233$ &$2.165$ \\
\cellcolor{gray!25}Reduced by (\%)  &\cellcolor{gray!25} -  &\cellcolor{gray!25}$0.05$ &\cellcolor{gray!25}$0.64$ &\cellcolor{gray!25}$0.64$  &\cellcolor{gray!25}$3.85$  &\cellcolor{gray!25}$0.69$ &\cellcolor{gray!25}$3.61$ &\cellcolor{gray!25}$1.00$ &\cellcolor{gray!25}$4.02$\\

\noalign{\hrule height .1pt} 
\multirow{2}{*}{}Min. Voltage (p.u.) 
       &$0.951$ &$0.955$ &$0.958$ &$0.953$ &$0.954$  &$0.967$ &$0.969$ &$0.975$ &$0.966$ \\
\cellcolor{gray!25}Increased by (p.u.)&\cellcolor{gray!25} -  &\cellcolor{gray!25} $0.004$ &\cellcolor{gray!25} $0.007$ &\cellcolor{gray!25} $0.002$ &\cellcolor{gray!25} $0.003$ &\cellcolor{gray!25}$0.016$ &\cellcolor{gray!25}$0.018$ &\cellcolor{gray!25}$0.024$ &\cellcolor{gray!25}$0.015$\\

\noalign{\hrule height .1pt}

\multirow{2}{*}{}Max. Voltage (p.u.) 
   &$1$ &$1.031$ &$1.001$ &$1.002$ &$1.001$  &$1.033$ &$1.027$ &$1.033$ &$1.007$ \\
\cellcolor{gray!25}Decreased by (p.u.)&\cellcolor{gray!25} -  &\cellcolor{gray!25} $-0.031$ &\cellcolor{gray!25} $-0.001$ &\cellcolor{gray!25} $-0.002$ &\cellcolor{gray!25} $-0.001$ &\cellcolor{gray!25}$-0.033$ &\cellcolor{gray!25}$-0.027$ &\cellcolor{gray!25}$-0.033$ &\cellcolor{gray!25}$-0.007$\\

\noalign{\hrule height .1pt}
\multirow{2}{*}{}Voltage Spread (p.u.) 
       &$0.036$ &$0.047$ &$0.030$ &$0.032$ &$0.031$  &$0.046$ &$0.040$ &$0.041$ &$0.023$ \\
\cellcolor{gray!25}Decreased by (p.u.)&\cellcolor{gray!25} -  &\cellcolor{gray!25} $-0.011$ &\cellcolor{gray!25} $0.006$ &\cellcolor{gray!25} $0.004$ &\cellcolor{gray!25} $0.005$ &\cellcolor{gray!25}$-0.010$ &\cellcolor{gray!25}$-0.004$ &\cellcolor{gray!25}$-0.005$ &\cellcolor{gray!25}$0.013$\\

\noalign{\hrule height .1pt}

  \end{tabular} \label{Default33}
\end{table*}

Table I presents the performance of the optimal VVO solution obtained under various case studies. In particular, $9$ cases, including a default test case, are simulated. Here, the default test case corresponds to the test system without any VVO equipment. The Table also presents the differences between the optimized system compared with the default test system in terms of the five metrics of quality. A negative sign in the row of comparisons indicates that optimizing the objective function (i.e., expected active power losses) in the system worsens that particular metric of quality.

Table I suggests that optimal capacitor switching and feeder reconfiguration can result in maximal improvements in the active power loss of the system. The results also show that optimal control of DSS marginally improves on the system losses as it decreases the peak load in the system. While feeder reconfiguration and optimal capacitor control can reduce the peak load of the system, the greatest reduction in the peak load comes from optimal DSS scheduling. In terms of minimum voltage of the feeder, the table shows that feeder reconfiguration in conjunction with ULTC can play a very important role. This impact on increasing the minimum voltage of the feeder becomes more significant if shunt capacitors are also employed for VVO in the system. 
Although the improvement in the minimum voltage reduces if all the equipment in the system are jointly controlled for power loss minimization, the average maximum voltage and the voltage spread improve significantly in return. Particularly, if reconfiguration switches, the ULTC, and shunt capacitors are jointly controlled in the system, the voltage spread in the feeder is $0.041$ p.u.. If the storage units are also employed for VVO, the voltage spread reduces to $0.013$ p.u. which is in agreement with a much better maximum voltage in this case.

Fig. 2 illustrates the voltage profile of the system under peak load of a typical day. Two test cases have been removed from the figure to make other curves more visible. Clearly, use of any VVO equipment in a day-ahead fashion has resulted in an improved voltage profile. Also, as expected, the smoothest voltage profile corresponds to the case where all VVO equipment are 
jointly operated in an optimal fashion.

\begin{figure}[btp!] 
\center
\includegraphics[height=2.6in, width = 3.4in]{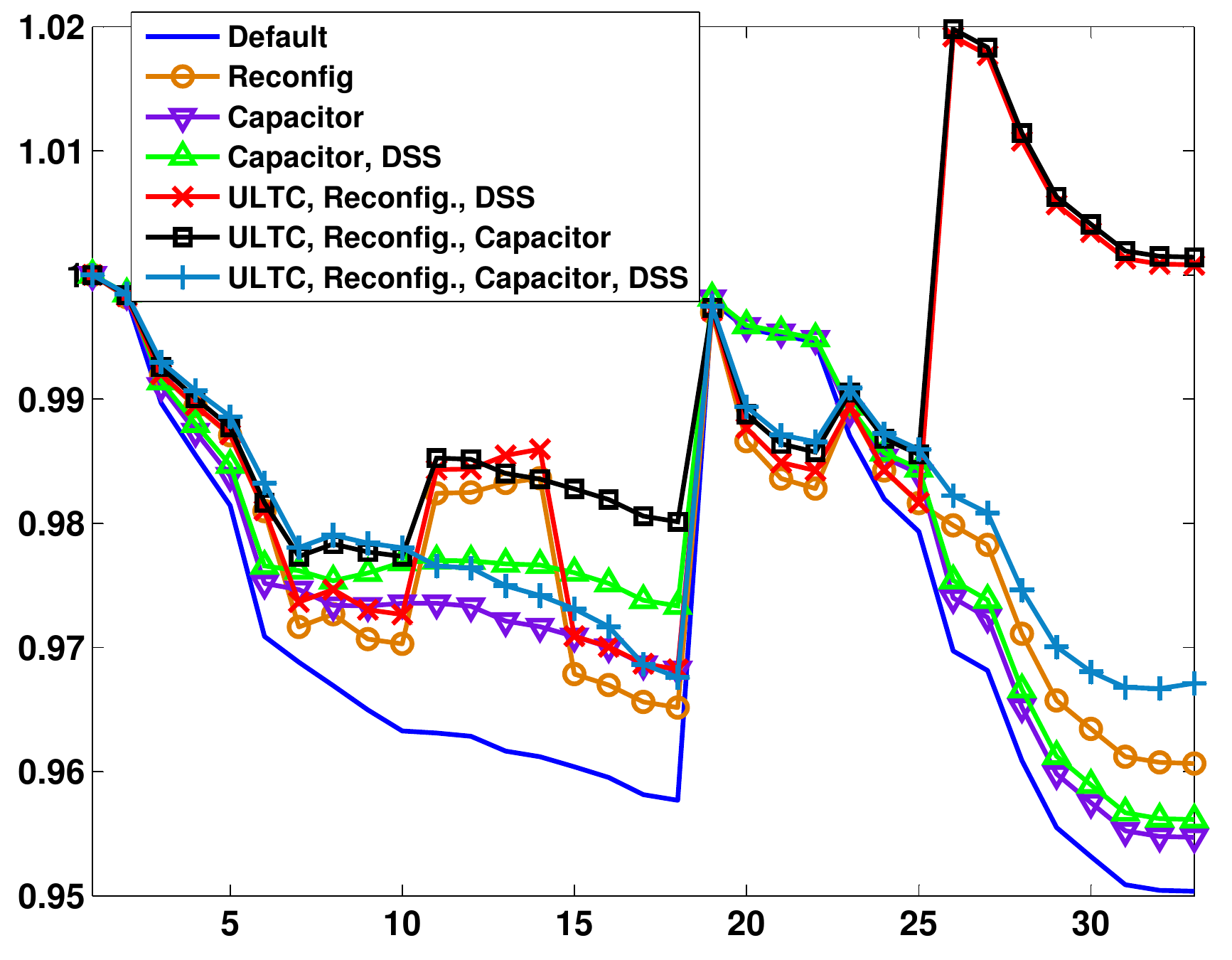}
\caption{Voltage profile of the system under different VVO test cases. The voltage profile is of a typical day at peak hour.}
\end{figure}

\section{Conclusion}\label{Conclusion}

In this paper, the optimal VVO problem aiming at minimizing the expected active power loss was formulated as a mixed-integer quadratic program to be solved using branch-and-bound methods. The methodology presented in the paper is comprehensive in that it jointly formulates the operation of shunt capacitors, storage units, reconfiguration switches, and ULTCs. The stochasticity of wind power generation was also addressed using a first-order Markov chain model. Numerical results on an active smart grid and using real data from smart meters and wind turbines was presented and discussed. The results indicated the significance of shunt capacitors and feeder reconfiguration in active power loss minimization and the importance of DSS in alleviating the peak load and improving the voltage profile of the system.


\ifCLASSOPTIONcaptionsoff
  \newpage
\fi

\bibliographystyle{IEEEtran}
\bibliography{IEEEabrv,ref}

\begin{thebibliography}{10}
\providecommand{\url}[1]{#1}
\csname url@samestyle\endcsname
\providecommand{\newblock}{\relax}
\providecommand{\bibinfo}[2]{#2}
\providecommand{\BIBentrySTDinterwordspacing}{\spaceskip=0pt\relax}
\providecommand{\BIBentryALTinterwordstretchfactor}{4}
\providecommand{\BIBentryALTinterwordspacing}{\spaceskip=\fontdimen2\font plus
\BIBentryALTinterwordstretchfactor\fontdimen3\font minus
  \fontdimen4\font\relax}
\providecommand{\BIBforeignlanguage}[2]{{%
\expandafter\ifx\csname l@#1\endcsname\relax
\typeout{** WARNING: IEEEtran.bst: No hyphenation pattern has been}%
\typeout{** loaded for the language `#1'. Using the pattern for}%
\typeout{** the default language instead.}%
\else
\language=\csname l@#1\endcsname
\fi
#2}}
\providecommand{\BIBdecl}{\relax}
\BIBdecl

\bibitem{High}
Y.~M. Atwa and E.~El-Saadany, ``Optimal allocation of {ESS} in distribution
  systems with a high penetration of wind energy,'' \emph{IEEE Trans. Power
  Syst.}, vol.~25, no.~4, pp. 1815--1822, 2010.

\bibitem{restoration}
C.~P. Nguyen and A.~J. Flueck, ``Agent based restoration with distributed
  energy storage support in smart grids,'' \emph{IEEE Trans. Smart Grid},
  vol.~3, no.~2, pp. 1029--1038, 2012.

\bibitem{LoadShaping}
T.~Jiang, Y.~Cao, L.~Yu, and Z.~Wang, ``Load shaping strategy based on energy
  storage and dynamic pricing in smart grid,'' \emph{IEEE Trans. Smart Grid},
  vol.~5, no.~6, pp. 2868--2876, 2014.

\bibitem{PrimaryVoltage}
K.~Christakou, D.-C. Tomozei, M.~Bahramipanah, J.-Y. Le~Boudec, and M.~Paolone,
  ``Primary voltage control in active distribution networks via broadcast
  signals: the case of distributed storage,'' \emph{IEEE Trans. Smart Grid},
  vol.~5, no.~5, pp. 2314--2325, 2014.

\bibitem{IntegratingElectrical}
P.~Wang, D.~H. Liang, J.~Yi, P.~F. Lyons, P.~J. Davison, and P.~C. Taylor,
  ``Integrating electrical energy storage into coordinated voltage control
  schemes for distribution networks,'' \emph{IEEE Trans. Smart Grid}, vol.~5,
  no.~2, pp. 1018--1032, 2014.

\bibitem{TheRole}
B.~P. Roberts and C.~Sandberg, ``The role of energy storage in development of
  smart grids,'' \emph{Proc. IEEE}, vol.~99, no.~6, pp. 1139--1144, 2011.

\bibitem{RobustMeter}
M.~Ghasemi~Damavandi, V.~Krishnamurthy, and J.~R. Mart\'i, ``Robust meter
  placement for state estimation in active distribution systems,'' \emph{IEEE
  Trans. Smart Grid}, vol.~6, no.~4, pp. 1972--1982, July 2015.

\bibitem{Valueof}
M.~R. Dorostkar-Ghamsari, M.~Fotuhi-Firuzabad, M.~Lehtonen, and A.~Safdarian,
  ``Value of distribution network reconfiguration in presence of renewable
  energy resources,'' \emph{IEEE Trans. on Power Syst.}, vol.~31, no.~3, pp.
  1879--1888, 2016.

\bibitem{SegmentedTime}
S.~Chen, W.~Hu, and Z.~Chen, ``Comprehensive cost minimization in distribution
  networks using segmented-time feeder reconfiguration and reactive power
  control of distributed generators,'' \emph{IEEE Trans. on Power Syst.},
  vol.~31, no.~2, pp. 983--993, 2016.

\bibitem{MHCjournal}
F.~Capitanescu, L.~F. Ochoa, H.~Margossian, and N.~D. Hatziargyriou,
  ``Assessing the potential of network reconfiguration to improve distributed
  generation hosting capacity in active distribution systems,'' \emph{IEEE
  Trans. Power Syst.}, vol.~30, no.~1, pp. 346--356, 2015.

\bibitem{SwitchingFrequency}
Z.~Li, S.~Jazebi, and F.~De~Leon, ``Determination of the optimal switching
  frequency for distribution system reconfiguration,'' \emph{IEEE Trans. on
  Power Delivery, to appear}.

\bibitem{TimeVarying}
R.~P. Broadwater, A.~H. Khan, H.~E. Shaalan, and R.~E. Lee, ``Time varying load
  analysis to reduce distribution losses through reconfiguration,'' \emph{IEEE
  Trans. on Power Delivery}, vol.~8, no.~1, pp. 294--300, 1993.

\bibitem{NetworkReconfiguration}
M.~E. Baran and F.~F. Wu, ``Network reconfiguration in distribution systems for
  loss reduction and load balancing,'' \emph{IEEE Trans. Power Del.}, vol.~4,
  no.~2, pp. 1401--1407, 1989.

\bibitem{reliability}
A.~Kavousi-Fard and T.~Niknam, ``Optimal distribution feeder reconfiguration
  for reliability improvement considering uncertainty,'' \emph{IEEE Trans.
  Power Del.}, vol.~29, no.~3, pp. 1344--1353, 2014.

\bibitem{restoration_reconfig}
A.~Botea, J.~Rintanen, and D.~Banerjee, ``Optimal reconfiguration for supply
  restoration with informed {A*} search,'' \emph{IEEE Trans. Smart Grid},
  vol.~3, no.~2, pp. 583--593, 2012.

\bibitem{ConvexModels}
J.~A. Taylor and F.~S. Hover, ``Convex models of distribution system
  reconfiguration,'' \emph{IEEE Trans. Power Syst.}, vol.~27, no.~3, pp.
  1407--1413, 2012.

\bibitem{PowerLoss}
R.~S. Rao, K.~Ravindra, K.~Satish, and S.~Narasimham, ``Power loss minimization
  in distribution system using network reconfiguration in the presence of
  distributed generation,'' \emph{IEEE Trans. Power Syst.}, vol.~28, no.~1, pp.
  317--325, 2013.

\bibitem{HamedAhmadi}
H.~Ahmadi, J.~R. Mart{\'\i}, and H.~W. Dommel, ``A framework for volt-{VAR}
  optimization in distribution systems,'' \emph{IEEE Trans. on Smart Grid},
  vol.~6, no.~3, pp. 1473--1483, 2015.

\bibitem{PeakLoadRelief}
A.~Padilha-Feltrin, D.~A.~Q. Rodezno, and J.~R.~S. Mantovani, ``{Volt-VAR}
  multiobjective optimization to peak-load relief and energy efficiency in
  distribution networks,'' \emph{IEEE Trans. on Power Delivery}, vol.~30,
  no.~2, pp. 618--626, 2015.

\bibitem{SSandVVO}
S.~Deshmukh, B.~Natarajan, and A.~Pahwa, ``State estimation and {voltage/VAR}
  control in distribution network with intermittent measurements,'' \emph{IEEE
  Trans. on Smart Grid}, vol.~5, no.~1, pp. 200--209, 2014.

\bibitem{SensitivityBased}
R.~A. Jabr and I.~D{\v{z}}afi{\'c}, ``Sensitivity-based discrete
  coordinate-descent for {volt/var} control in distribution networks,''
  \emph{IEEE Trans. on Power Syst., to appear}, 2016.

\bibitem{IntegrationDGVVO}
J.~Barr and R.~Majumder, ``Integration of distributed generation in the
  {Volt/VAR} management system for active distribution networks,'' \emph{IEEE
  Trans. on Smart Grid}, vol.~6, no.~2, pp. 576--586, 2015.

\bibitem{VVO1}
I.~Roytelman, B.~Wee, and R.~Lugtu, ``Volt/var control algorithm for modern
  distribution management system,'' \emph{IEEE Trans. on Power Syst.}, vol.~10,
  no.~3, pp. 1454--1460, 1995.

\bibitem{ThroughDG}
S.~Deshmukh, B.~Natarajan, and A.~Pahwa, ``{Voltage/VAR} control in
  distribution networks via reactive power injection through distributed
  generators,'' \emph{IEEE Trans. on Smart Grid}, vol.~3, no.~3, pp.
  1226--1234, 2012.

\bibitem{VVONiknam}
T.~Niknam, M.~Zare, and J.~Aghaei, ``Scenario-based multiobjective volt/var
  control in distribution networks including renewable energy sources,''
  \emph{IEEE Trans. on Power Delivery}, vol.~27, no.~4, pp. 2004--2019, 2012.

\bibitem{Acomprehensive}
F.~Capitanescu, I.~Bilibin, and E.~R. Ramos, ``A comprehensive centralized
  approach for voltage constraints management in active distribution grid,''
  \emph{IEEE Trans. on Power Syst.}, vol.~29, no.~2, pp. 933--942, 2014.

\bibitem{MISOCP_VAR}
Z.~Tian, W.~Wu, B.~Zhang, and A.~Bose, ``Mixed-integer second-order cone
  programing model for {VAR} optimisation and network reconfiguration in active
  distribution networks,'' \emph{IET Generation, Transmission \& Distribution},
  vol.~10, no.~8, pp. 1938--1946, 2016.

\bibitem{MCMC}
G.~Papaefthymiou and B.~Klockl, ``{MCMC} for wind power simulation,''
  \emph{IEEE Trans. on Energy Conversion}, vol.~23, no.~1, pp. 234--240, 2008.

\bibitem{CapacitorPlacement}
M.~E. Baran and F.~F. Wu, ``Optimal sizing of capacitors placed on a radial
  distribution system,'' \emph{IEEE Trans. Power Del.}, vol.~4, no.~1, pp.
  735--743, 1989.

\bibitem{Relief}
A.~Padilha-Feltrin, D.~A.~Q. Rodezno, and J.~R.~S. Mantovani, ``{Volt-VAR}
  multiobjective optimization to peak-load relief and energy efficiency in
  distribution networks,'' \emph{IEEE Trans. on Power Delivery}, vol.~30,
  no.~2, pp. 618--626, 2015.

\bibitem{Integration}
G.~Carpinelli, G.~Celli, S.~Mocci, F.~Mottola, F.~Pilo, and D.~Proto, ``Optimal
  integration of distributed energy storage devices in smart grids,''
  \emph{IEEE Trans. Smart Grid}, vol.~4, no.~2, pp. 985--995, 2013.

\bibitem{InverterLess}
Z.~Wang, H.~Chen, J.~Wang, and M.~Begovic, ``Inverter-less hybrid voltage/var
  control for distribution circuits with photovoltaic generators,'' \emph{IEEE
  Trans. on Smart Grid}, vol.~5, no.~6, pp. 2718--2728, 2014.

\bibitem{TimeInterval}
Z.~Hu, X.~Wang, H.~Chen, and G.~Taylor, ``{Volt/VAr} control in distribution
  systems using a time-interval based approach,'' \emph{IEE
  Proceedings-Generation, Transmission and Distribution}, vol. 150, no.~5, pp.
  548--554, 2003.

\bibitem{HybridVVO}
A.~Mohapatra, P.~R. Bijwe, and B.~K. Panigrahi, ``An efficient hybrid approach
  for volt/var control in distribution systems,'' \emph{IEEE Trans. on Power
  Delivery}, vol.~29, no.~4, pp. 1780--1788, 2014.

\bibitem{GridSupport}
M.~Nick, R.~Cherkaoui, and M.~Paolone, ``Optimal allocation of dispersed energy
  storage systems in active distribution networks for energy balance and grid
  support,'' \emph{IEEE Trans. Power Syst.}, vol.~29, no.~5, pp. 2300--2310,
  2014.

\bibitem{FlexibleOptimal}
A.~Gabash and P.~Li, ``Flexible optimal operation of battery storage systems
  for energy supply networks,'' \emph{IEEE Trans. Power Syst.}, vol.~28, no.~3,
  pp. 2788--2797, 2013.

\bibitem{EnergyStorageDevices}
L.~H. Macedo, J.~F. Franco, M.~J. Rider, and R.~Romero, ``Optimal operation of
  distribution networks considering energy storage devices,'' \emph{IEEE Trans.
  on Smart Grid}, vol.~6, no.~6, pp. 2825--2836, 2015.

\bibitem{Imposing}
M.~Lavorato, J.~F. Franco, M.~J. Rider, and R.~Romero, ``Imposing radiality
  constraints in distribution system optimization problems,'' \emph{IEEE Trans.
  Power Syst.}, vol.~27, no.~1, pp. 172--180, 2012.

\bibitem{node_32}
M.~E. Baran and F.~F. Wu, ``Network reconfiguration in distribution systems for
  loss reduction and load balancing,'' \emph{IEEE Trans. Power Del.}, vol.~4,
  no.~2, pp. 1401--1407, 1989.

\bibitem{CER}
\BIBentryALTinterwordspacing
{Commission for Energy Regulation (CER)}. {Customer Behaviour Trials (CBTs)}.
  [Online]. Available: \url{http://www.cer.ie/electricity-gas/smart-metering/}
\BIBentrySTDinterwordspacing

\bibitem{ISSDA}
\BIBentryALTinterwordspacing
{Irish Social Science Data Archive}. {CER} smart metering project. [Online].
  Available:
  \url{http://www.ucd.ie/issda/data/commissionforenergyregulationcer/}
\BIBentrySTDinterwordspacing

\bibitem{WindData}
\BIBentryALTinterwordspacing
{EirGrid}. [Online]. Available:
  \url{http://www.eirgrid.com/operations/systemperformancedata/}
\BIBentrySTDinterwordspacing

\bibitem{BES}
K.~Divya and J.~{\O}stergaard, ``Battery energy storage technology for power
  systemsóan overview,'' \emph{Elect. Power Syst. Res.}, vol.~79, no.~4, pp.
  511--520, 2009.

\bibitem{MATPOWER}
R.~D. Zimmerman, ``{MATPOWER}: A {MATLAB} power system simulation package.''

\bibitem{CVX}
M.~Grant, S.~Boyd, and Y.~Ye, ``{CVX}: Matlab software for disciplined convex
  programming,'' 2008.

\end{thebibliography}

\end{document}